\documentclass[10pt]{amsart}
\usepackage{url,graphicx}
\usepackage[dvips]{epsfig}
\usepackage{latexsym,color}
\usepackage{amscd,amssymb,verbatim}
\newcommand{\be}{\begin{enumerate}}
\newcommand{\ee}{\end{enumerate}}

\newcommand{\nin}{\noindent}
\newcommand{\pr}{\noindent{\bf Proof. }}
\newcommand{\sm}{\setminus}

\newcommand{\ca}{{\mathcal A}}
\newcommand{\cc}{{\mathcal C}}
\newcommand{\cf}{{\mathcal F}}
\newcommand{\cj}{{\mathcal J}}
\newcommand{\cl}{{\mathcal L}}
\newcommand{\cn}{{\mathcal N}}

\newcommand{\zz}{{{\mathbb Z}_2}}

\newcommand{\sea}{\searrow}
\newcommand{\nea}{\nearrow}

\newcommand{\bd}{{\text{\rm Bd}\,}}
\newcommand{\bl}{{\text{\rm Bl}}}

\newcommand{\da}{\Delta}
\newcommand{\dgn}{\text{\rm DG}_{\text{\rm n}}}

\newcommand{\fix}{\text{\rm Fix}\,}

\newcommand{\im}{\text{\rm Im}}
\newcommand{\lk}{\text{\rm lk}}
\newcommand{\lov}{{\mathcal Lo}}
\newcommand{\op}{\text{op}}
\newcommand{\ra}{\rightarrow}
\newcommand{\sd}{{\text{\rm sd}\,}}
\newcommand{\st}{\text{\rm st}}
\newcommand{\thom}{\text{\tt Hom}\,}
\newcommand{\ti}{\tilde}

\newtheorem{thm}{Theorem}[section]
\newtheorem{df}[thm]{Definition}

\newtheorem{crl}[thm]{Corollary}
\newtheorem{prop}[thm]{Proposition}
\newtheorem{conj}[thm]{Conjecture}
\newtheorem{rem}[thm]{Remark}

\numberwithin{equation}{section}
\numberwithin{figure}{section}

\begin{document}

\title[Simple homotopy types of some combinatorially defined complexes]
{Simple homotopy types of Hom-complexes, neighborhood complexes, 
Lov\'asz complexes, and atom crosscut complexes}

\author{Dmitry N. Kozlov}
\address{Institute of Theoretical Computer Science, Eidgen\"ossische Technische
Hochschule, Z\"urich, Switzerland}
\email{dkozlov@inf.ethz.ch}
\thanks {Research supported by Swiss National Science Foundation Grant PP002-102738/1}

\date\today

\begin{abstract}
  In this paper we provide concrete combinatorial formal deformation
  algorithms, namely sequences of elementary collapses and expansions,
  which relate various previously extensively studied families of
  combinatorially defined polyhedral complexes.

  To start with, we give a~sequence of elementary collapses leading
  from the barycentric subdivision of the neighborhood complex to the
  Lov\'asz complex of a~graph.  Then, for an~arbitrary lattice $\cl$
  we describe a~formal deformation of the barycentric subdivision of
  the atom crosscut complex $\Gamma(\cl)$ to its order
  complex~$\da(\bar\cl)$. We proceed by proving that the complex of
  sets bounded from below $\cj(\cl)$ can also be collapsed
  to~$\da(\bar\cl)$.

  Finally, as a~pinnacle of our project, we apply all these results to
  certain graph complexes. Namely, by describing an explicit formal
  deformation, we prove that, for any graph $G$, the neighborhood
  complex $\cn(G)$ and the polyhedral complex $\thom(K_2,G)$ have the
  same simple homotopy type in the sense of Whitehead.
\end{abstract}

\maketitle

\section{Introduction.}

Motivation for the research presented in this paper came from the
quest for better understanding of the relationship between
neighborhood complexes and \thom-complexes.

Originally, neighborhood complexes were introduced and used by
Lov\'asz, see \cite{Lo}, to attack the Kneser Conjecture, as well as
to provide some of the first nontrivial algebro-topological lower
bounds for chromatic numbers of graphs. After an~active period of
research and several attempts at the generalizations of the
neighborhood complexes, the so-called \thom-complexes were introduced,
again by Lov\'asz. We refer the reader to the survey
article~\cite{IAS}.

\thom-complexes depend on two parameters, both of them graphs. One of 
the motivations for introducing these gadgets was the fact that the
polyhedral complex $\thom(K_2,G)$ turned out to be homotopy equivalent
to the simplicial complex $\cn(G)$, for any graph $G$, see e.g.,
\cite[Proposition 4.2]{BK03b} for an~argument. We remark, that all
known proofs of this fact are in a~way nonconstructive, making use of
statements like Quillen's Fiber Lemma.

One of the classical views of topology is combinatorial, using
``moves'' between cell complexes called elementary collapses and
elementary expansions, see e.g., \cite{Al} for a~prototypical
approach. These ideas were further developed and reached their
maturity in the work of Whitehead, see e.g.,~\cite{Wh}.

The suggested modus operandi would be to, instead of looking for
continuous homotopies, construct a~discrete object: a~so-called formal
deformation.  The natural question of whether two homotopy equivalent
spaces would necessarily be connected by a~formal deformation turned
out to have a~negative answer, and as a~result an~exciting and
important theory of simple homotopy type and Whitehead torsion ensued.

Precisely this circle of ideas has been the driving force behind this
article. As a~consequence, we were able to find explicit formal
deformations between various combinatorially defined simplicial
complexes. These sequences, when concatenated, prove that
$\thom(K_2,G)$ and $\cn(G)$ have the same simple homotopy type, for
any graph~$G$.

In the process of constructing these formal deformations we had to
revise and upgrade several central results from Topological
Combinatorics, these are Theorems \ref{thm:cnlov}
and~\ref{thm:crossc}. Classically these results would just conclude
the existence of a~deformation retraction, or, even sometimes only the
existence of a~homotopy equivalence, see~\cite{Bj}. Here, both in
Theorem \ref{thm:cnlov}, and in Theorem \ref{thm:crossc}, we provide
an~explicitly described algorithmic sequence of collapses, and
expansions (expansions are only required in the second theorem).

Furthermore, we needed to consider a~complex, which does not seem to
have appeared before. For any finite lattice we construct a~simplicial
complex of bounded below sets of elements: its vertices are elements
of $\bar\cl$ and $S\subseteq\bar\cl$ is a~simplex if and only if $S$
has lower bound different from~$\hat 0$. It turned out that this
complex collapses onto the order complex $\da(\bar\cl)$, and that this
sequence of collapses can be described algorithmically.

Finally, in the last section we bring all these results into play in
order to construct the promised formal deformation between the studied
graph complexes. It is well known, see \cite{Co}, that in a~formal
deformation all the expansions can be carried out first, followed by
the collapses, however we have chosen to present our formal
deformation as in~\eqref{eq:conc}, since it passes through several
complexes, which appear to be of interest in their own.

Unfortunately, our formal deformation is still rather complicated.
Finding a~simpler natural formal deformation from $\thom(K_2,G)$ to
$\cn(G)$ remains a~challenging task.

\section{Notations.}

We start by recalling some notations. For a partially ordered set $P$
we let $\da(P)$ denote its {\it order complex} (also known as the {\it
nerve} of the corresponding category), that is the simplicial complex
whose set of vertices is the set of elements of $P$, and whose set of
simplices is the set of {\it chains} (ordered subsets) of~$P$.

Let $P$ be an arbitrary partially ordered set. For any subset $S\subseteq P$, we let $P[S]$ denote the induced partial order on $S$. We let $P^\op$ denote the poset whose set of elements is the same as that of $P$, but whose partial order is the reverse of the partial order of~$P$. Note that for arbitrary poset $P$ we have $\da(P)=\da(P^\op)$. The minimal, resp.\ maximal, element of $P$ (if it has one), is denoted by $\hat 0$, resp.\ $\hat 1$. In this case, we set $\bar P:=P\sm\{\hat 0, \hat 1\}$. If $P$ has a~minimal element, we let $\ca(P)$ denote the set {\it atoms} of $P$, i.e., of elements which cover $\hat 0$.

An order-preserving map $\varphi:P\ra P$ (i.e.\ a~map such that $x\geq
y$ implies $\varphi(x)\geq\varphi(y)$) is called a~{\it monotone map},
if for any $x\in P$ either $x\geq\varphi(x)$ or $x\leq\varphi(x)$. If
$x\geq\varphi(x)$ for all $x\in P$, then we call $\varphi$ a~{\it
decreasing map}, analogously, if $x\leq\varphi(x)$ for all $x\in P$,
then we call $\varphi$ an~{\it increasing map}.

For an arbitrary lattice $\cl$, and a~subset $S\subseteq\cl$, we let $\bigwedge S$ denote the common meet of all the elements in $S$, and, analogously, we let $\bigvee S$ denote the common join of all the elements in $S$.

\section{Barycentric and stellar subdivisions.}

\label{subd}

For an~arbitrary CW complex $X$, we let $\cf(X)$ denote its {\it face
poset}: the partially ordered set whose elements are all nonempty
cells of $X$, and whose partial order is given by the cell
inclusion. When $X$ is a~regular CW complex, we let $\bd(X)$ denote
its {\it barycentric subdivision}. Clearly, we have
$\bd(X)=\da(\cf(X))$.

For a~simplicial complex $X$, and an~arbitrary simplex $\sigma\in X$,
let $\lk_X\sigma$ denote the {\it link} of $\sigma$ in $X$, let
$\st_X\sigma$ denote the {\it closed star} of~$\sigma$ in $X$.
Furthermore, let $\sd(X,\sigma)$ denote the {\it stellar subdivision}
of $X$ at $\sigma$. The effect that the stellar subdivision has on the
face poset $\cf(X)$ is a~special case of the combinatorial blowup in
a~lattice: $\cf(\sd(X,\sigma))\cup\{\hat
0\}=\bl_\sigma(\cf(X)\cup\{\hat 0\})$,
see~\cite[Proposition~4.9]{FK1}.

It is a~classical fact that the barycentric subdivision can be
represented as a~sequence of stellar subdivisions: simply take
a~reverse linear extension of $\cf(X)$ and perform stellar
subdivisions of the corresponding simplices in this
order. Combinatorially, using the terminology of \cite{FK1}, this
corresponds to taking the whole poset $\cf(X)$ as
a~building set.  The above mentioned fact can then be seen as
a~special case of
\cite[Theorem 3.4]{FK1}.

When $Y$ is a~simplicial subcomplex of $X$, we say that $X$ {\it
collapses onto} $Y$ if there exists a~sequence of elementary collapses
leading from $X$ to $Y$; in this case we write $X\searrow Y$ (or,
equivalently, $Y\nearrow X$); we refer the reader to \cite[\S 4, p.\
14]{Co}, for the definition of the elementary collapse for an
arbitrary finite CW pair. The reverse of an elementary collapse is
called an elementary {\it expansion}. A~sequence of elementary
collapses and elementary expansions leading from a~complex $X$ to the
complex $Y$ is called a~{\it formal deformation}. If such a~sequence
exists, then the simplicial complexes $X$ and $Y$ are said to have the
same {\it simple homotopy type}, see~\cite{Al,Co,Wh}.

It is well-known, see e.g., \cite[\S 25, Statement (25.1)]{Co}, that
a~subdivision of any CW complex $X$ has the same simple homotopy type
as $X$. For completeness of our results we describe here an
explicit formal deformation from $X$ to $\bd X$.

To start with, since the barycentric subdivision can be represented as
a~sequence of stellar subdivisions, it is enough to find a formal
deformation leading from $X$ to $\sd(X,\sigma)$, for an arbitrary
simplex $\sigma\in X$. One choice of such deformation is
a~concatenation of two steps.

\vskip3pt

\nin{\it Deformation algorithm from $X$ to $\sd(X,\sigma)$.}

\vskip3pt

\nin {\bf Step 1.} Add a~cone over $\st_X(\sigma)$. More precisely, 
consider a new simplicial complex $X'$, such that $V(X')=V(X)\cup\{v\}$, 
$X$ is an~induced subcomplex of $X'$, and $\lk_{X'}v=\st_X(\sigma)$.

\vskip3pt

\nin {\bf Step 2.} Delete from $X'$ all the simplices containing $\sigma$.

\vskip3pt

Since $\st_X(\sigma)$ is a~cone, in particular collapsible, the Step~1
can be performed as a~sequence of elementary expansions. Furthermore,
the Step~2 can be performed as a~sequence of elementary collapses as
follows. The set of the simplices which are to be deleted can be
written as a~disjoint union of sets $A$ and $B$, where $B$ is the set
of all simplices which contain both $\sigma$ and $v$. Clearly, adding
$v$ to a~simplex is a~bijection $\mu:A\ra B$. Let
$\{\tau_1\dots,\tau_t\}$ be a~reverse linear extension order on $A$,
then $\{(\tau_1,\mu(\tau_1),\dots,(\tau_t,\mu(\tau_t))\}$ is
an~elementary collapsing sequence.

Finally, we see that performing Steps 1 and 2, in this order, will
yield a~stellar subdivision of $X$ at $\sigma$, and therefore our
description is completed.

\section{Collapsing the neighborhood complex of a~graph onto its Lov\'asz complex.}

The next theorem is a~specialization of \cite[Theorem 3.1(b)]{Ko1} to
the case of the finite posets.

\begin{thm} \label{thm1} \cite{Ko1}.
Let $P$ be a finite~poset, and let $\varphi:P\ra P$ be a~monotone map. Assume $P\supseteq Q\supseteq\fix\varphi$, then $\Delta(P)$ collapses onto $\Delta(Q)$.
\end{thm}

Let $G$ be an~arbitrary undirected graph. We let $V(G)$ denote the set
of vertices of $G$. For any $v\in V(G)$, we let $N(v)$ denote the set
of all neighbors of $v$, i.e., $N(v)=\{w\in V(G)\,|\,(w,v)\in
E(G)\}$. Furthermore, for any subset $S\subseteq V(G)$, we set
$N(S):=\bigcap_{v\in S} N(v)$, i.e., $N(S)$ denotes the set of common
neighbors of all the vertices in~$S$.

In \cite{Lo}, Lov\'asz has introduced the following class of simplicial complexes, in order to study the topological obstructions to graph colorings.

\begin{df} \label{df:nbhd}
For an arbitrary graph $G$, let $\cn(G)$ be the simplicial complex, whose set of vertices consists of all nonisolated vertices of $G$, and whose set of simplices consists of all subsets $S\subseteq V(G)$, such that the vertices in $S$ have a~common neighbor, i.e., such that $N(S)\neq\emptyset$.
\end{df}

These complexes have been studied fairly extensively, see e.g.,~\cite{Cs,Ziv}.

Note that $N$ induces an~order-reversing map
$N:\cf(\cn(G))\ra\cf(\cn(G))$, in particular $N^2(A)\supseteq A$, for
any $A\subseteq V(G)$. It can also be seen that $N^3=N$.

\begin{figure}[hbt]
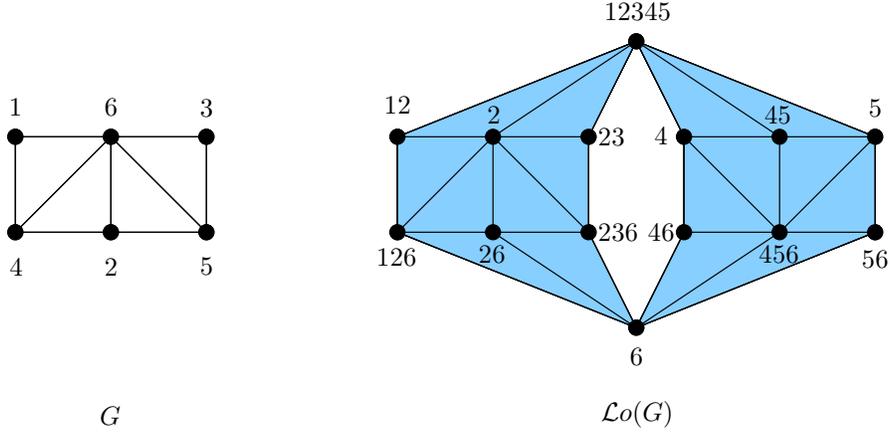

\begin{center}
  \begin{picture}(0,0)%
    \includegraphics{ex1.pstex}%
  \end{picture}%
  \input{ex1.pstex_t}%
  
\end{center}
\caption{A graph $G$ and its Lov\'asz complex.}
\label{fig:ex1}
\end{figure}

\begin{df}
For an arbitrary graph $G$, the complex $\da(N(\cf(\cn(G))))$ is
called the {\bf Lov\'asz complex} of $G$ and is denoted by $\lov(G)$.
\end{df}

One property, which distinguishes the Lov\'asz complex as an
interesting object of study, is that it possesses a~natural
$\zz$-action, induced by the map $N$.  Indeed, we see that
$$N(\cf(\cn(G)))\supseteq N^2(\cf(\cn(G)))\supseteq
N^3(\cf(\cn(G)))=N(\cf(\cn(G))),$$ hence
$N(\cf(\cn(G)))=N^2(\cf(\cn(G)))$.  It is an easy check that when the
graph~$G$ has no loops, this action is free, and so in this case
$\lov(G)$ has a~natural structure of $\zz$-space.

It is well-known that for any graph, its neighborhood complex and its
Lov\'asz complex are homotopy equivalent. The next proposition
strengthens this result.

\begin{thm}\label{thm:cnlov}
The simplicial complex $\bd(\cn(G))$ collapses onto the simplicial complex $\lov(G)$.
In particular, $\cn(G)$ and $\lov(G)$ have the same simple homotopy type.
\end{thm}

\pr
Define the map $\varphi:\cf(\cn(G))\ra\cf(\cn(G))$, by simply setting
$\varphi:=N^2$. From our previous comments, it is clear that $\varphi$
is an order-preserving map, and that $A\leq\varphi(A)$, for any
$A\subseteq V(G)$. Note, that it is also true that
$\varphi^2=\varphi$, but we do not need this additional fact.

We conclude that $\varphi$ is an~ascending map, and hence, by the
Theorem~\ref{thm1} we obtain that the simplicial complex
$\da(\cf(\cn(G)))=\bd(\cn(G))$ collapses onto the simplicial complex
$\da(\varphi(\cf(\cn(G))))=\da(N(\cf(\cn(G))))=\lov(G)$.
\qed


\begin{rem} 
By the discussion in Section~\ref{subd}, the Theorem \ref{thm:cnlov}
allows us to construct an~explicit formal deformation from $\cn(G)$
to~$\lov(G)$.
\end{rem}

\section{Simple homotopy type of crosscut complexes.}

Crosscut complexes play a prominent role in Topological Combinatorics, e.g., see the survey \cite{Bj}.

\begin{df} \label{df:crossc}
Let $\cl$ be a lattice, the {\bf atom crosscut complex} $\Gamma(\cl)$ associated to $\cl$ is a~simplicial complex defined as follows:
\begin{itemize}
\item the set of vertices of $\Gamma(\cl)$ is equal to the set of atoms of $\cl$, in other words, $V(\Gamma(\cl))=\ca(\cl)$;
\item the subset $\sigma\subseteq \ca(\cl)$ is a simplex in $\Gamma(\cl)$ if and only if
the join of elements in $\sigma$ is not equal to~$\hat 1$.
\end{itemize}
\end{df}

Recall that a lattice $\cl$ is called {\it atomic}, if all elements of $\cl$ can be represented as joins of atoms. For an arbitrary lattice $\cl$, let $\cl_a$ denote the sublattice
consisting of $\hat 0$, and of all the elements which are joins of atoms.

\begin{thm}\label{thm:crossc}
Let $\cl$ be an arbitrary finite lattice.
\begin{enumerate}
\item [(a)] If $\cl$ is atomic, then the simplicial complex $\bd(\Gamma(\cl))$ collapses onto the simplicial complex $\da(\bar\cl)$.
\item [(b)] In the general case, both $\bd(\Gamma(\cl))$ and $\da(\bar\cl)$ collapse onto the simplicial complex~$\da(\bar\cl_a)$.
\end{enumerate}
In both cases we conclude that the simplicial complexes $\Gamma(\cl)$ and $\da(\bar\cl)$ have the same simple homotopy type.
\end{thm}

\pr
Assume first that $\cl$ is atomic. Define a map
$\varphi:\cf(\Gamma(\cl))\ra\cf(\Gamma(\cl))$ as follows: a~simplex
$\sigma$ is mapped to $\ca(\cl)_{\leq \bigvee \sigma}$. To start with,
the map $\varphi$ is well-defined, since
$\bigvee\ca(\cl)_{\leq\bigvee\sigma}\leq\bigvee\sigma<\hat 1$, also,
clearly $\varphi(\sigma)\supseteq\sigma$. Furthermore, $\varphi$ is
order-preserving, since if $\tau\supseteq\sigma$, then
$\bigvee\tau\geq\bigvee\sigma$, implying
$\ca(\cl)_{\leq\bigvee\tau}\supseteq\ca(\cl)_{\leq\bigvee\sigma}$.  We
remark, that
$\bigvee\sigma\geq\bigvee\varphi(\sigma)\geq\bigvee\sigma$, hence
$\bigvee\sigma=\bigvee\varphi(\sigma)$, and therefore
$\varphi^2(\sigma)=\varphi(\sigma)$; however we do not need the latter
fact for our argument.

From the discussion above we see that $\varphi$ is a~monotone map, and hence, by the Theorem \ref{thm1} we conclude that the simplicial complex $\da(\cf(\Gamma(\cl)))=\bd(\Gamma(\cl))$ collapses onto the simplicial complex $\da(\varphi(\cf(\Gamma(\cl))))$. On the other hand, since the lattice is atomic, we have $\varphi(\cf(\Gamma(\cl)))=\bar\cl$, and so, as desired, the simplicial complex $\bd(\Gamma(\cl))$ collapses onto the simplicial complex $\da(\bar\cl)$.
 
Now, remove the assumption that $\cl$ is atomic, and consider the general case. By the argument above we see that $\bd(\Gamma(\cl))$ collapses onto $\da(\bar\cl_a)$. On the other hand, it is not difficult to check that the order-preserving map $\psi:\cl\ra\cl$ mapping $x$ to the join of the elements of $\ca(\cl)_{\leq x}$, and mapping $\hat 0$ to itself, is a~descending map. Its image is precisely~$\cl_a$.
\qed


\begin{rem} 
Again, by the discussion in Section~\ref{subd}, the
Theorem~\ref{thm:crossc} can be used to construct an~explicit formal
deformation from $\Gamma(\cl)$ to $\da(\bar\cl)$.
\end{rem}

Due to its general nature, the Theorem \ref{thm:crossc} has many applications. Let us mention one of them.

\begin{df}
Let $n$ be any natural number, and let $\dgn$ be the simplicial complex of all disconnected graphs on $n$ labeled vertices.  In other words, the vertices of $\dgn$ are all
pairs $(i,j)$, with $i<j$, $i,j\in [n]$, i.e., all possible edges of
a~graph on labeled $n$ vertices; and simplices of $\dgn$ are all
collections of edges which form a graph with at least 2 connected
components.
\end{df}

Recall, that for an arbitrary natural number $n$, $\Pi_n$ denotes the {\it partition lattice}: the poset consisting of all set partitions of the set $\{1,\dots,n\}$, partially ordered by partition refinement.

\begin{crl}
The simplicial complex $\bd(\dgn)$ collapses onto $\da(\bar\Pi_n)$.
\end{crl}

\pr A direct check yields $\dgn=\Gamma(\Pi_n)$, hence the result follows from the Theorem~\ref{thm:crossc}.
\qed

\vskip5pt

We remark that the complex $\dgn$ appeared in the work of Vassiliev on
knot theory, \cite{Va}, whereas $\da(\bar\Pi_n)$ encodes the geometry
of the braid arrangement by means of the Goresky-MacPherson theorem,
see~\cite{GM}.

\vskip5pt

Recall, that for an arbitrary lattice $\cl$, a~{\it crosscut} is
a~subset $C\subseteq\bar\cl$, such that:
\begin{itemize} 
\item $C$ is an antichain (a set of mutually incomparable elements);
\item $C$ is saturated in the following sense: for any chain $\gamma$ 
of $\cl$ there exists an element $x\in C$, such that $\gamma\cup\{x\}$
is again a~chain.
\end{itemize}

Generalizing the Definition~\ref{df:crossc} the {\it crosscut complex}
$\Gamma(C,\cl)$ associated to the crosscut $C$ is a~simplicial complex
defined as follows:
\begin{itemize}
\item the set of vertices of $\Gamma(C,\cl)$ is equal to the set $C$;
\item the subset $\sigma\subseteq C$ is a simplex in $\Gamma(C,\cl)$ 
if and only if either the join of the elements in $\sigma$ is not
equal to $\hat 1$, or the meet of the elements in $\sigma$ is not
equal to~$\hat 0$.
\end{itemize}

The set of atoms is a~special case of a crosscut, and the atom
crosscut complex is a~special case of the crosscut complex.

Naturally, a crosscut $C$ divides the lattice $\cl$ into two parts
$$\cl_{\geq C}=\{x\in\cl\,|\,x\geq s,\text{ for some }s\in C\},$$ and
$$\cl_{\leq C}=\{x\in\cl\,|\,x\leq s,\text{ for some }s\in C\},$$
which intersect in $C$. Let $\cl_C$ be the subposet consisting of
$\hat 0$, $\hat 1$, and of all joins and meets of the elements of the
crosscut $C$. Let $\varphi:\cl\ra\cl$ be a~map defined as follows:
\[\varphi(x)=
\begin{cases}
\bigvee C_{\leq x},&\text{ if }x\in\cl_{\geq C};\\
\bigwedge C_{\geq x},&\text{ if }x\in\cl_{\leq C}.
\end{cases}\]
We can see that $\varphi$ is order-preserving. The only nontrivial
case to be checked is when $x\geq y$, $x\in\cl_{\geq C}$, and
$y\in\cl_{\leq C}$. Since in this case $x\geq y$ is a chain, there
must exist an element $z\in C$, such that $\{x,y,z\}$ is also
a~chain. Obviously, we must have $x\geq z\geq y$. Since $z\in C_{\leq
x}\cap C_{\geq y}$, we conclude that $z\geq\varphi(y)$, and
$\varphi(x)\geq z$, hence $\varphi(x)\geq\varphi(y)$.

It is also easy to check that $\varphi$ is a~monotone map, namely
$\varphi(x)\leq x$, if $x\in\cl_{\geq C}$, and $\varphi(x)\geq x$, if
$x\in\cl_{\leq C}$. Furthermore, the image of $\varphi$ is precisely
$\cl_C$. By Theorem~\ref{thm1} we see that $\da(\bar\cl)$ collapses
onto $\da(\bar\cl_C)$.

Interestingly, Sonja \v{C}uki\'c has remarked that in general the
simplicial complex $\bd(\Gamma(C,\cl))$ does not have to collapse onto
the simplicial complex~$\da(\bar\cl_C)$, \cite{Cu1}. We conclude this
section by conjecturing that the weak version of
Theorem~\ref{thm:crossc} is still true in general.

\begin{conj}
For an arbitrary lattice $\cl$ and an arbitrary crosscut $C$, the
simplicial complex $\bd(\Gamma(C,\cl))$ and the simplicial
complex~$\da(\bar\cl)$ have the same simple homotopy type.
\end{conj}

Together with our previous observations, this conjecture can equivalently be formulated as:

\begin{conj}
For an arbitrary lattice $\cl$ and an arbitrary crosscut $C$, the
simplicial complex $\Gamma(C,\cl)$ and the simplicial
complex~$\da(\bar\cl_C)$ have the same simple homotopy type.
\end{conj}


\section{Collapsing the complex of sets bounded from below onto the order complex.}

We start by defining a~combinatorial gadget, which provides
a~convenient language for describing sequences of elementary
collapses.

\begin{df}
Let $P$ be a poset with the covering relation $\succ$.
\begin{itemize} 
\item We define a~{\bf~partial matching} on $P$ to be a~set 
$\Sigma\subseteq P$, and an~injective map $\mu:\Sigma\ra P\sm\Sigma$,
such that $\mu(x)\succ x$, for all $x\in\Sigma$. 
\item The elements of $P\sm(\Sigma\cup\mu(\Sigma))$ are called {\bf critical}.
We let $\cc(P,\mu)$ denote the set of critical elements.
\item Additionally, such a partial matching $\mu$ is called {\bf acyclic} if there
exists no sequence of distinct elements $x_1,\dots,x_t\in\Sigma$,
where $t\geq 2$, satisfying $\mu(x_1)\succ x_2$, $\mu(x_2)\succ x_3$,
$\dots$, $\mu(x_t)\succ x_1$.
\end{itemize}
\end{df}

The partial acyclic matchings and elementary collapses are closely
related, as the next proposition shows.

\begin{prop}\label{DMT}
  Let $\da$ be a~regular CW complex and $\da'$ a~subcomplex of $\da$,
  then the following are equivalent:

  a) there is a~sequence of elementary collapses leading from $\da$
  to~$\da'$;

  b) there is a~partial acyclic matching on the poset $\cf(\da)$ with
  the set of critical cells being exactly $\cf(\da')$.
\end{prop}
\pr See \cite[Proposition 5.4]{Koz2}.
\qed\vspace{5pt}

We remark that the implication b)$\Rightarrow$a) is a~special case of
a~more general result proved by R.\ Forman, see~\cite{Fo}.

There is a number of constructions associating a~simplicial complex to
a~poset (or more generally, to a~category), here is one which works
for lattices.

\begin{df} \label{df:jcomp}
Let $\cl$ be an~arbitrary finite lattice. We define $\cj(\cl)$ be the
simplicial complex whose set of vertices is equal to the set of
elements of $\bar\cl$, and whose simplices are all subsets
$S\subseteq\bar\cl$ which have a~nontrivial lower bound, i.e., such
that $\bigwedge S\neq\hat 0$.
\end{df}

Clearly, the simplicial complex $\cj(\cl)$ contains $\da(\bar\cl)$ as
a~subcomplex.  It turns out that much more is true.

\begin{thm} \label{thm:main}
Let $\cl$ be an~arbitrary finite lattice, then $\cj(\cl)\sea\da(\bar\cl)$.
\end{thm}

\pr As the centerpiece of the argument we define the following partial 
acyclic matching on $\cf(\cj(\cl))$. Let $S$ be an arbitrary simplex
of $\cj(\cl)$. Assume that $\cf(\cj(\cl))[S]$ is not a~chain. Set
$t:=|S|$, and let $\{a_1,a_2,\dots,a_t\}$ be a~linear extension of
$\cf(\cj(\cl))[S]$, i.e., if $1\leq i<j\leq t$, then $a_i\not\ge a_j$.

\begin{figure}[hbt]
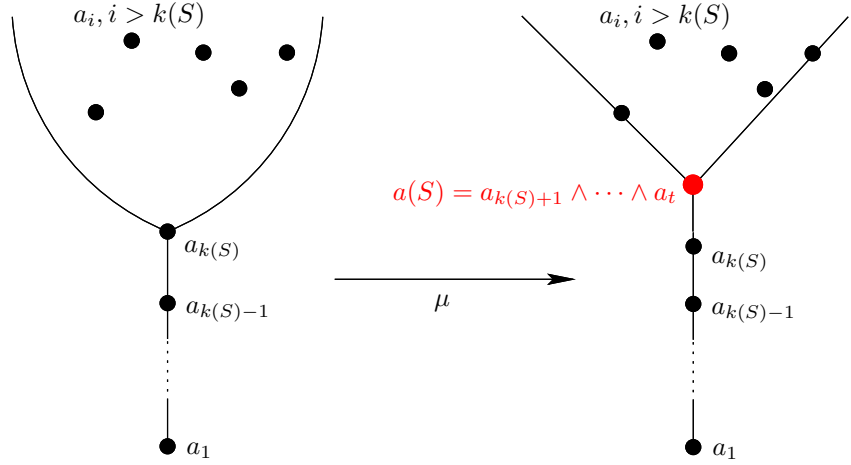

\begin{center}
  \begin{picture}(0,0)%
    \includegraphics{mu.pstex}%
  \end{picture}%
  \input{mu.pstex_t}%
  
\end{center}
\caption{The partial matching $\mu$.}
\label{fig:mu}
\end{figure}

Let $k(S)$ be the maximal index, $1\leq k(S)\leq t$, such that
$a_1<a_2<\dots<a_{k(S)}$, and $a_{k(S)}<a_i$, for all $k(S)+1\leq
i\leq t$, see Figure~\ref{fig:mu}. If $S$ has no minimal element, then
we set $k(S):=0$. Set $a(S):=a_{k(S)+1}\wedge\dots\wedge a_t$. Since
$\cf(\cj(\cl))[S]$ is not a~chain, we have $k(S)\leq t-2$, and hence
$a(S)$ is well-defined.

Let $\Sigma$ be the set of all subsets $S\subseteq\bar\cl$, such that
$\cf(\cj(\cl))[S]$ is not a~chain, and such that $a(S)\not\in S$. For
$S\in\Sigma$ define $\mu(S):=S\cup\{a(S)\}$, again see
Figure~\ref{fig:mu}. Clearly, $\mu$ defines a~partial matching, and,
since for any $S\in\Sigma$ we have $a(\mu(S))=a(S)$, we see that the
set $\mu(\Sigma)\cup\Sigma$ consists of all subsets
$S\subseteq\bar\cl$, such that $\cf(\cj(\cl))[S]$ is not
a~chain. Consequently, the set of critical elements
$\cc(\cf(\cj(\cl)),\mu)$ consists of all chains
$S\in\cf(\da(\bar\cl))$.

Let us see that the partial matching $\mu$ is acyclic. Assume there
exists a~sequence $S_1,\dots,S_t\in\Sigma$, where $t\geq 2$, such that
$\mu(S_1)\succ S_2$, $\mu(S_2)\succ S_3$, $\dots$, $\mu(S_t)\succ
S_1$. Let again $\{a_1,a_2,\dots,a_t\}$ be a~linear extension of
$\cf(\cj(\cl))[S_1]$, as above. By the definition of covering
relations, and, since $S_2\neq S_1$, we have $S_2=\mu(S_1)\sm\{a_i\}$,
for some $1\leq i\leq t$.  If $1\leq i\leq k(S_1)$, then
$a(S_2)=a(S_1)$, which, together with $S_1=\mu(S_1)\sm\{a(S_1)\}$,
implies $a(S_2)\in S_2$, and hence $S_2\in\mu(\Sigma)$, giving
a~contradiction.

Finally, the only option left is that $k(S_1)+1\leq i\leq t$, in which
case $a(S_2)\geq a(S_1)$, since the join is taken over a~set, where
each element is larger than $a(S_1)$. If the equality $a(S_2)=a(S_1)$
holds, then $S_2\in\mu(\Sigma)$, again giving a~contradiction. Thus we
have shown that a~strict inequality must hold: $a(S_2)> a(S_1)$.

Analogously, we can prove that $a(S_{i+1})> a(S_i)$, for all $1\leq
i\leq t-1$, and that $a(S_1)> a(S_t)$, which, when combined together,
yields a~contradiction to the assumption that the matching is not
acyclic. By Proposition \ref{DMT} we see that the acyclic matching
$\mu$ provides a~sequence of elementary collapses leading from
$\cj(\cl)$ to~$\da(\bar\cl)$.
\qed

\vskip5pt

We invite the interested reader to see what the statement of the
Theorem~\ref{thm:main} translates to for their favorite lattice~$\cl$.

\section{Application to graph complexes.}

Let $G$ and $T$ be two undirected graphs. Recall that the set map
$\varphi:V(G)\ra V(T)$ is called a~{\it graph homomorphism} from $G$
to $T$ if, for any pair of vertices $x,y\in V(G)$, such that $(x,y)\in
E(G)$, we have $(\varphi(x),\varphi(y))\in E(T)$.

\begin{df} \label{dfhom}
For arbitrary undirected graphs $T$ and $G$, we let $\thom(T,G)$
denote the~polyhedral complex whose cells are indexed by all functions
$\eta:V(T)\rightarrow 2^{V(G)}\setminus\{\emptyset\}$, such that for
any $(x,y)\in E(T)$, we have $\eta(x)\times\eta(y)\subseteq E(G)$.

The closure of a~cell $\eta$ consists of all cells indexed by
functions $\ti\eta:V(T)\rightarrow 2^{V(G)}\setminus\{\emptyset\}$,
which satisfy $\ti\eta(v)\subseteq\eta(v)$, for all $v\in V(T)$.
\end{df}

We note that the set of vertices of $\thom(T,G)$ coincides with the
set of all graph homomorphisms from $T$ to $G$, so the polyhedral
complex $\thom(T,G)$ may be thought of as an~appropriate
topologization of this set.

The $\thom$-complexes were introduced by Lov\'asz, and recently
studied in a~series of papers, see
\cite{BK03a,BK03b,BK03c,CK1,CK2,Ko2,IAS,Ziv}, in connection with topological
obstructions to graph colorings.

For the case $T=K_2$, the Definition \ref{dfhom} can be restated
somewhat more directly.  Recall that, for arbitrary $A,B\subseteq
V(G)$, $A,B\neq\emptyset$, we call the pair $(A,B)$ a~{\it complete
bipartite subgraph} of $G$, if for any $x\in A$, $y\in B$, we have
$(x,y)\in E(G)$, i.e., $A\times B\subseteq E(G)$. Let $\Delta^{V(G)}$
be the simplex whose set of vertices is $V(G)$, in particular, the
faces of $\Delta^{V(G)}$ can be identified with the subsets of~$V(G)$.

Clearly, $\Delta^{V(G)}\times\Delta^{V(G)}$ is a~polyhedral complex,
whose cells are direct products of two simplices.  $\thom(K_2,G)$ can
be identified as the subcomplex of $\Delta^{V(G)}\times\Delta^{V(G)}$
defined by the following condition: $\sigma\times\tau\in\thom(K_2,G)$
if and only if $(\sigma,\tau)$ is a~complete bipartite subgraph
of~$G$.

We are now ready to formulate one of the main results of this paper.

\begin{thm} \label{thm:appl}
For an arbitrary graph $G$, the neighborhood complex $\cn(G)$ and the
polyhedral complex $\thom(K_2,G)$ have the same simple homotopy type.
\end{thm}

\pr Set $P:=\cf^\op(\thom(K_2,G))\cup\{\hat 0,\hat 1\}$. 
As was mentioned before, $P$ is a~lattice, and $\da(\bar P)
=\bd\thom(K_2,G)$.  By the Theorem~\ref{thm:crossc}(b), we see that
both simplicial complexes $\bd\thom(K_2,G)$ and $\bd\Gamma(P)$
collapse onto the simplicial complex $\da(\bar P_a)$.

\begin{figure}[hbt]
\begin{center}
  \begin{picture}(0,0)%
    \includegraphics{ex2.pstex}%
  \end{picture}%
  \input{ex2.pstex_t}%
  
\end{center}
\caption{The poset $\bar P_a$, for $P=\cf^\op(\thom(K_2,G))\cup\{\hat 0,\hat 1\}$.}
\label{fig:ex2}
\end{figure}

\vskip3pt

\nin {\it Description of $\Gamma(P)$}. The vertices of $\Gamma(P)$ are all 
the pairs $(A,B)$, $A,B\subseteq V(G)$, such that $N(A)=B$, and
$N(B)=A$. These can be indexed with the simplices $A\in\cn(G)$,
$A\in\im N$, which is the same as to take the elements of
$N(\cf(\cn(G)))$, or the vertices of $\da(N(\cf(\cn(G))))=\lov(G)$.

The simplices of $\Gamma(P)$ are all sets of pairs
$\{(A_1,B_1),\dots,(A_t,B_t)\}$, such that $\bigcap_{i=1}^t
A_i\neq\emptyset$, and $\bigcap_{i=1}^t B_i\neq\emptyset$. Since
$N(A)\cap N(B)=N(A\cup B)$, for arbitrary subsets $A,B\subseteq V(G)$,
and since $B_i=N(A_i)$, for $1\leq i\leq t$, the second condition
amounts to saying that $N(\bigcup_{i=1}^t A_i)\neq\emptyset$.

\vskip3pt

Let $\cl$ denote the poset of all $A\in\cn(G)$, $A\in\im N$, ordered
by inclusion, with a~minimal and a~maximal elements attached. Clearly,
$\da(\bar\cl)=\lov(G)$. From the description of $\Gamma(P)$ above, we
see that $\da(\bar\cl)$ is a~subcomplex of $\Gamma(P)$. On the other
hand, by the Theorem~\ref{thm:main}, the simplicial complex $\cj(\cl)$
collapses onto $\da(\bar\cl)$.

Let $\mu$ be the acyclic matching from the proof of the
Theorem~\ref{thm:main} which gives the collapsing sequence. We claim
that the restriction of $\mu$ to $\cf(\Gamma(P))$ is again an~acyclic
matching. Since $\cf(\Gamma(P))$ is a~lower ideal in $\cf(\cj(\cl))$,
the only thing which has to be checked is that if
$S\in\cf(\Gamma(P))\cap\Sigma$, then $\mu(S)\in\cf(\Gamma(P))$; here
$\Sigma$ is as in the proof of the Theorem~\ref{thm:main}.

Assume that $S=\{A_1,\dots,A_t\}$, where the sets are listed in the
linear extension order, i.e., if $1\leq i<j\leq t$, then
$A_i\not\supseteq A_j$.  Let $a(S)$ be the subset of $V(G)$ defined as
in the proof of the Theorem~\ref{thm:main}. Clearly, $a(S)\subseteq
A_t$, this implies that $a(S)\cup\bigcup_{i=1}^t A_i=\bigcup_{i=1}^t
A_i$, and therefore, the set of pairs
$$\mu(S)=\{(A_1,B_1),\dots,(A_t,B_t),(a(S),N(a(S)))\}$$ is a~simplex of
$\Gamma(P)$.

We conclude that the restriction of $\mu$ to $\cf(\Gamma(P))$ gives
a~collapsing sequence from $\Gamma(P)$ to $\lov(G)$.

Let us summarize our findings in the following concatenation of
sequences of collapses and expansions:
\begin{equation}\label{eq:conc}
\bd\thom(K_2,G)\sea\da(\bar P_a)\nea\bd\Gamma(P),\quad
\Gamma(P)\sea\lov(G)\nea\bd\cn(G),
\end{equation}
where the first two sequences are given by the
Theorem~\ref{thm:crossc}(b), the third sequence is given by the restriction
of the acyclic matching $\mu$ as above, and the fourth sequence is given by
the Theorem~\ref{thm:cnlov}.

The discussion in Section \ref{subd} implies now that the polyhedral
complex of all bipartite subgraphs of $G$, $\thom(K_2,G)$, and the
neighborhood complex $\cn(G)$, have the same simple homotopy type, and yields
an~explicit formal deformation between these two complexes.
\qed


\begin{rem}
All 4 sequences of collapses and expansions can be nondegenerate.
Figures~\ref{fig:ex1} and~\ref{fig:ex2} show an example of a~graph which satisfies this.
\end{rem}

\vskip10pt

\nin {\bf Acknowledgments.} We would like to thank Sonja \v{C}uki\'c for 
the careful reading of the initial draft of this paper, and for
helping to improve the presentation of our results. We are indebted to
the Swiss National Science Foundation and ETH-Z\"urich for the
financial support of this research.

\end{document}